\documentclass{article}

\usepackage{amsthm}
\usepackage{amsmath}
\usepackage{amssymb}
\usepackage{graphicx}
\usepackage{epsfig}

\newtheorem{theorem}{Theorem}

\newtheorem{definition}[theorem]{Definition}
\newtheorem{example}[theorem]{Example}
\newtheorem{remark}[theorem]{Remark}
\newtheorem{corollary}[theorem]{Corollary}

\newcommand{\examp}[1]
  {\begin{example} {\rm #1} \end{example}}

\begin{document}

\begin{center} \LARGE{\textbf{DUALLY DEGENERATE VARIETIES AND THE GENERALIZATION OF A  THEOREM OF GRIFFITHS--HARRIS}}

\vspace*{3mm}

{\large  Maks A. Akivis and Vladislav V. Goldberg}

\end{center}

Abstract. {\footnotesize The dual variety $X^*$ for a smooth
$n$-dimensional variety $X$ of the projective space $\mathbb{P}^N$
is the set of tangent hyperplanes to $X$. In the general case, the
variety $X^*$ is a hypersurface in the dual space
$(\mathbb{P}^{N})^*$. If $\dim X^* < N - 1$, then the variety $X$
is called dually degenerate.

The authors refine these definitions for a variety $X \subset
\mathbb{P}^{N}$ with a degenerate Gauss map of rank $r$. For such
a variety, in the general case, the dimension of its dual variety
$X^*$ is $N - l - 1$, where $l = n - r$, and $X$ is dually
degenerate if $\dim X^* < N - l - 1$.

In 1979 Griffiths and Harris proved that a smooth variety $X
\subset \mathbb{P}^{N}$ is dually degenerate if and only if all
its second fundamental forms are singular. The authors generalize
this theorem for  a variety $X \subset \mathbb{P}^{N}$ with a
degenerate Gauss map of rank $r$.}

\vspace*{2mm}


\vspace*{2mm}


\vspace*{3mm}

\setcounter{equation}{0}

\section{Dual varieties}

\textbf{1.1 Varieties with Degenerate Gauss Maps.} An almost
everywhere smooth $n$-dimensional variety $X$ of a projective
space $\mathbb{P}^N$ is called  {\em tangentially degenerate} or
{\em a variety with a degenerate Gauss map} if the rank of its
Gauss map
$$
\gamma: X \rightarrow \mathbb{G} (n, N)
$$
is less than $n, \; 0 \leq r = \mbox{{\rm rank}} \; \gamma < n$.
Here $n = \dim X, \, x \in X, \, \gamma (x) = T_x (X)$, and $T_x
(X)$ is the tangent subspace to $X$ at $x$ considered as an
$n$-dimensional projective space $\mathbb{P}^n$.  The number  $r$
is also called the {\em rank} of $X, \; r = \mbox{{\rm rank}} \;
X$. The case $r = 0 $ is trivial one: it gives just an $n$-plane.
The varieties with degenerate Gauss maps were studied in many
books and papers (see, for example, [AG 02], [AG 04], [FP 01], [GH
79], [L 99], where one can find further references).

Let $X \subset  \mathbb{P}^N$ be  an $n$-dimensional almost
everywhere smooth variety with a degenerate Gauss map. Suppose
that  $0 < \mbox{{\rm rank}} \, \gamma = r < n$. Denote by $L$ a
leaf of this map, $L = \gamma^{-1} (T_x) \subset X;\, \dim L = n -
r = l$. The number $l = n - r$ is called the \textit{Gauss defect}
of the variety $X$. The leaves $L$ of the map $\gamma$ are
$l$-dimensional subspaces of the  space $\mathbb{P}^N, \, L = P^l
\subset \mathbb{P}^N$ or open parts of such subspaces (see, for
example, [AG 04], Theorem 3.1, p. 95, or Theorem 2.10 in [GH~79];
or \S 5 in [L 99]; or the Linearity Theorem in Section 2.3 of [FP
01]).

A variety with a degenerate Gauss map of rank $r$ foliates into
its leaves  $L$ of dimension $l$, along which the tangent subspace
$T_x (X)$ is fixed.

The tangent subspace $T_x (X)$ is fixed when a point $x$ moves
along  regular points of $L$. This is the reason that we denote it
by $T_L, \, L \subset T_L$. A pair $(L, T_L)$ on $X$ depends on
$r$ parameters.

The foliation on $X$ defined as indicated above is called the {\em
Monge--Amp\`{e}re foliation}.

The varieties of rank $r < n$ are multidimensional analogues of
developable surfaces of a three-dimensional Euclidean space.

The main  results on the geometry of varieties with degenerate
Gauss maps and further references can be found in Chapter 4 of the
book [AG 93] and in the recently published book [AG 04].

\textbf{1.2 Dual Defect and Dually Degenerate Varieties.} By the
duality principle, to a point $x$ of a projective space
$\mathbb{P}^N$, there corresponds a hyperplane $\xi$. The set of
hyperplanes of the space $\mathbb{P}^N$ forms the dual projective
space $(\mathbb{P}^{N})^*$ of the same dimension $N$. Under this
correspondence, to a subspace $\mathbb{P} \subset \mathbb{P}^N$ of
dimension $p$, there corresponds a subspace $\mathbb{P}^* \subset
(\mathbb{P}^{N})^*$ of dimension $N - p - 1$. Under the dual map,
the incidence of subspaces is reversed, that is, if $\mathbb{P}_1
\subset \mathbb{P}_2$, then $\mathbb{P}_1^* \supset
\mathbb{P}_2^*$.

Let $X$ be  an irreducible, almost everywhere smooth variety of
dimension $n$  in the space $\mathbb{P}^N$, let $x$ be a smooth
point of $X$, and let $T_x X$ be the tangent subspace to $X$ at
the point $x$. A hyperplane $\xi$ is said to be \textit{tangent to
$X$ at $x$} if $T_x \subset \xi$. The bundle of hyperplanes $\xi$
tangent to $X$ at $x$ is of dimension $N - n - 1$.

The set of all hyperplanes $\xi$ tangent to the variety $X$ at its
smooth points composes a variety
$$
X^{\wedge} = \{\xi \subset \mathbb{P}^N | \exists \, x \in X_{sm}
\;\, \text{such that} \;\, T_x X \subseteq \xi\},
$$
where $X_{sm}$ is the locus of smooth points of the variety $X$.
But this variety can be not closed if $X$ has singular points. The
\textit{dual variety} $X^*$ of a variety $X$ is the closure of the
variety $X^{\wedge}$:
\begin{equation}\label{eq:1}
X^* = \overline{X^{\wedge}} = \overline{\{\xi \subset \mathbb{P}^N
| \exists \, x \in X_{sm} \;\, \text{such that} \;\, T_x X
\subseteq \xi\}}.
\end{equation}

The dual variety $X^*$ can also be described as the envelope of
the family of hyperplanes $\xi$ dual to the points $x \in X$. This
gives a practical way for finding $X^*$, which we will use in
examples.

If a variety $X$ is tangentially nondegenerate, i.e., if its rank
$r = n$, then in the general case, the dimension $n^*$ of its dual
variety $X^*$ is equal to
\begin{equation}\label{eq:2}
n^* =  \dim X^* = (N - n - 1) + n = N - 1.
\end{equation}
Equation (2) means that the variety $X^*$ is a hypersurface with a
degenerate Gauss map in the space $(\mathbb{P}^{N})^*$. The rank
$r$ of $X^*$ equals the dimension $n$ of the variety $X, \, r =
\text{rank} \, X^* = n$, and its Gauss defect
$$
\delta_\gamma (X^*)
= l^* = n^* - r = N - r - 1.
$$

However, it may happen that $\dim X^* < N - 1$. Then the number
\begin{equation}\label{eq:3}
\delta_* = N - 1 - \dim X^*
\end{equation}
is called the \textit{dual defect} of the variety $X$, and the
variety $X$ itself is said to be  \textit{dually degenerate}. The
classification of dually degenerate smooth varieties of small
dimensions $n$ with positive dual defect $\delta_*$ was found in
[E 85, 86] for $n \leq 6$; in [E 85, 86] and [LS 87] for $n = 7$,
and in [BFS 92] for $n \leq 10$ (see also Section 9.2.C in
[T~01]).

The dual defect of a variety $X$ must be defined as the difference
between an expected dimension of the dual variety $X^*$ and its
true dimension. An expected dimension of the dual variety $X^*$ of
smooth (tangentially nondegenerate) varieties equals $N - 1$. For
these reasons, the above standard definitions of the dual defect
and dually degenerate varieties are appropriate for smooth
varieties.

However, in the books [FP 01] (p. 55); [Ha 92] (p. 199); [L 99]
(p. 16); and [T~01], the above standard definitions of the dual
defect and dually degenerate varieties, which are appropriate for
tangentially nondegenerate varieties, are automatically extended
to varieties with degenerate Gauss maps. In our opinion, these
definitions should be refined, because for these latter varieties,
the expected dimension of $X^*$ is less than $N -1$, and for them
the appropriate definition of the dual defect (and dually
degenerate varieties) must be different (see below).

If a variety $X$ has a degenerate Gauss map (i.e., if its rank $r
< n$), then the dual variety $X^*$ is a fibration whose fiber is
the bundle $\Xi = \{\xi \subset \mathbb{P}^N | \xi \supseteq T_L
X\}$ of hyperplanes $\xi$ containing the tangent subspace $T_L X$
and whose base is the manifold $B = X^* / \Xi$. The dimension of a
fiber $\Xi$ of this fibration (as in the case $r = n$) equals $N -
n - 1, \, \dim \Xi = N - n - 1$, and the dimension of the base $B$
equals $r, \, \dim B = r$, i.e., the dimension of $B$ coincides
with the rank of the variety $X$. Therefore, in the general case,
the dimension $n^*$ of its dual variety $X^*$ is determined by the
formula
\begin{equation}\label{eq:4}
 \dim X^* = (N - n - 1) + r = N - l - 1,
\end{equation}
where $l = \dim L = \delta_\gamma = n - r$, and its Gauss defect
is equal to $\delta_\gamma (X^*) = l^* \linebreak = n^* - r =  (N
- l - 1) - r = N - n - 1 = \dim \Xi.$

Note that formula (4) for an expected dimension of the dual
variety of a variety with degenerate Gauss map appeared also in
the paper [P 02]  and implicitly in the books [L 99] (see 7.2.1.1
and 7.3i) and [FP 01] (see Section 2.3.4).

However, it may happen that $\dim X^* < N - l - 1$. Then the
number
\begin{equation}\label{eq:5}
\delta_* = N - l - 1 - \dim X^*
\end{equation}
is called the \textit{refined dual defect} of the variety $X$, and
the variety $X$ itself is said to be  \textit{dually degenerate}.

We will make the following three remarks about these  definitions
of the refined dual defect and dually degenerate varieties with
degenerate Gauss maps:

\begin{itemize}

\item[(i)] The refined dual defect of tangentially nondegenerate
varieties can be obtained from this  definition (5) by taking $l =
0$.

\item[(ii)] While by the standard definition (3), all varieties
with degenerate Gauss maps are dually degenerate, by the
definition (5), they can be either dually degenerate or dually
nondegenerate. Moreover, while by the standard definition (3), the
dual defect $\delta_*$ of a dually nondegenerate variety with
degenerate Gauss map equals its Gauss defect, $\delta_* =
\delta_\gamma = n - r > 0$, by the definition (5), the refined
dual defect $\delta_*$ of such a variety equals $0$, $\delta_* =
0$, and this is more appropriate for a \textit{dually
nondegenerate} variety.

\item[(iii)] Note also  that dually degenerate smooth varieties in
the projective space $\mathbb{P}^N$ are few and far between. As to
dually degenerate varieties with degenerate Gauss maps, we are
aware of only a few examples of  dually degenerate varieties $X$
with  degenerate Gauss maps: the varieties $X$ with  degenerate
Gauss maps of ranks three and four in $\mathbb{P}^N$ were
considered in [P 02]. At the end of this paper we construct
another example of a dually degenerate variety with a degenerate
Gauss map (see Example 9).
\end{itemize}

The following theorem follows immediately from the preceding
considerations.

\begin{theorem}
Let $X$ be a dually nondegenerate variety with a degenerate Gauss
map of dimension $n$ and rank $r$ in the space $\mathbb{P}^N$.
Then the leaves $L$ of the Monge--Amp\`{e}re foliation of $X$ are
of dimension $l = n - r$. The dual variety $X^* \subset
(\mathbb{P}^{N})^*$ is of dimension
\begin{equation}\label{eq:6}
n^* = N - l - 1
\end{equation}
and the same rank $r$, and the leaves $L^*$ of the
Monge--Amp\`{e}re foliation of $X^*$ are of dimension
\begin{equation}\label{eq:7}
l^* = N - n - 1.
\end{equation}
Under the Gauss map, the plane generator $L^*$ corresponds to a
tangent subspace $T_x (X)$ of the variety $X$, and  the tangent
subspace $T_\xi (X^*)$ of the variety $X^*$ corresponds to a plane
generator $L$, i.e., on $X$ the tangent bundle $T (X)$ and the
Monge--Amp\`{e}re foliation $L (X)$ are mutually dual.
\end{theorem}

In particular, if a variety $X \subset \mathbb{P}^N$ is
tangentially nondegenerate, then we have $n = r, \, l = 0$ (i.e.,
$n^* = N - 1$), and the dual map (*) sends $X$ to a hypersurface
$X^* \subset (\mathbb{P}^{N})^*$ with a degenerate Gauss map of
rank $n$ with the leaves $L^*$ of the Monge--Amp\`{e}re foliation
of dimension $l^* = N - n - 1$.

Conversely, if $X$ is a hypersurface with a degenerate Gauss map
of rank $r < N - 1$ in $\mathbb{P}^N$, then the variety $X^*$ dual
to $X$ is a tangentially nondegenerate variety of dimension $r$
and rank $r$.

In particular, the dual map (*) sends a tangentially nondegenerate
variety $X \subset \mathbb{P}^N$ of dimension and rank $r = n = N
- 2$ to a hypersurface $X^* \subset (\mathbb{P}^{n+2})^*$ with a
degenerate Gauss map of rank $r$, and $X^*$ bears an $r$-parameter
family of rectilinear generators. Each of these rectilinear
generators possesses $r$ foci if each is counted as many times as
its multiplicity. The  hypersurface $X^*$ is torsal and foliates
into $r$ families of torses (for definition of the torse, see
Example 5 in Section 1.3). The original variety $X$ bears a net of
conjugate lines corresponding to the torses of the variety $X^*$.
Of course, the correspondence indicated above is mutual.

We consider an irreducible, almost everywhere smooth variety $X$
of dimension $n$ and rank $r$ in the space $\mathbb{P}^N$ in more
detail. The tangent bundle $T (X)$ of $X$ is formed by the
$n$-dimensional subspaces $T_x$ tangent to $X$ at points $x \in X$
and depending on $r$ parameters. The subspaces $T_x$ are tangent
to $X$ along the plane generators $L$ of dimension $l = n - r$
composing on $X$ the Monge--Amp\`{e}re foliation $L (X)$. The
bundle $T (X)$ and the foliation $L (X)$ have a common
$r$-dimensional base.

Let (*) be the dual map of $\mathbb{P}^N$ onto $(\mathbb{P}^N)^*$.
The dual map (*) sends the variety $X$ to a variety $X^*$, which
is the set of all hyperplanes $\xi \subset (\mathbb{P}^N)^*$
tangent to $X$ along the leaves $L$  of its Monge--Amp\`{e}re
foliation. The map (*) sends the tangent bundle $T (X)$ and the
Monge--Amp\`{e}re foliation $L (X)$ of $X$  to the
Monge--Amp\`{e}re foliation $L (X^*)$ and the tangent bundle $T
(X^*)$ of $X^*$, respectively. Thus, under the dual map (*), we
have
$$
(T (X))^* = L (X^*), \;\;\; (L (X))^* = T (X^*),
$$
where $\dim T (X^*) = \dim X^* = n^* = N - l - 1$ and $\dim L
(X^*) = \dim L^* = l^* \linebreak = N - n - 1$.

\textbf{1.3 Examples.} We now consider a few examples. In
particular, we consider the main types of varieties with
degenerate Gauss maps (cones, multidimensional torses, and
tangentially degenerate hypersurfaces) and determine their dual
varieties. Most of these examples can be found in [AG 02] and [AG
04], Section 2.4. We present them here in order to illustrate the
notion of the dual variety.

\setcounter{theorem}{0}

\examp{\label{examp:1} First, we consider a simple example. Let
$X$ be a smooth spatial curve $X$ in a three-dimensional
projective space $\mathbb{P}^3$. For this curve, we have $N = 3,
\, n = r = 1, \, l = 0$, and  $T_x (X)$ is the tangent line to $X$
at $x$. The dual map (*) sends a point $x \in X$ to a plane $\xi
\subset X^*$, and the dual variety $X^*$ is the envelope of the
one-parameter family of hyperplanes $\xi$, i.e., $X^*$ is a torse.

Using formulas (6) and (7) for $n^*$ and $l^*$,  we find that $n^*
= 2, \, l^* = 1$. The variety $X^*$ bears rectilinear generators
$L^*$ along which the tangent planes $\xi = T(X^*)$ are constant.
Hence $\text{rank} \; X^* = 1$. The generators $L^*$ of the torse
$X^*$ are dual to the tangent lines $T (X)$ to the curve $X$. }

\examp{\label{examp:2} Let $X$ be a hypersurface with a degenerate
Gauss map of rank $r < N - 1$ in $\mathbb{P}^N$. In this case, we
have $l = N - 1 - r, \, n^* = (N - 1) - l = (N - 1) - (N - 1 - r)
= r, \, l^* = n^* - r = 0$, i.e., the variety $X^*$ dual to $X$ is
a tangentially nondegenerate variety of dimension $r$ and rank
$r$.

Conversely, if a variety $X \subset \mathbb{P}^N$ is tangentially
nondegenerate, then we have $n = r, \, l = 0$, and the dual map
(*) sends $X$ to a hypersurface $X^* \subset (\mathbb{P}^{N})^*$
(i.e., $n^* = N - 1$) with a degenerate Gauss map of rank $n$ with
the leaves $L^*$ of the Monge--Amp\`{e}re foliation of dimension
$l^* = N - n - 1$.}

\examp{\label{examp:3}  In the space $\mathbb{P}^N, \; N \geq 4$,
we take two arbitrary smooth space curves, $Y_1$ and $Y_2$, that
do not belong to the same three-dimensional space, and the set of
all straight lines intersecting these two curves. These straight
lines form a three-dimensional variety $X$. Such a variety is
called the \emph{join}. Its dimension is three, $n = \dim \, X =
3$. It is easy to see that the  variety  $X$ has a degenerate
Gauss map. In fact, the three-dimensional tangent subspace $T_x
(X)$ to $X$ at a point $x$ lying on a rectilinear generator $L$ is
defined by this generator $L$ and two straight lines tangent to
the curves $Y_1$ and $Y_2$ at the points $y_1$ and $y_2$ of their
intersection with the line $L$. Because this tangent subspace does
not depend on the location of the point $x$ on the  generator $L$,
the  variety  under consideration is a variety  $X = V_2^3$ with a
degenerate Gauss map of rank two. Thus, in this case, we have $l =
1, \, r = r^* = 2, \, n^* = N - 2, \, l^* = N - 4$.}

This example can be generalized by taking $k$ space curves in the
space $\mathbb{P}^N$, where $N \geq 2k$  and $k > 2$, and
considering a $k$-parameter family of $(k-1)$-planes intersecting
all these $k$ curves.

\examp{\label{examp:4}  Suppose that $S$ is a subspace of the
space $\mathbb{P}^N, \; \dim S = l - 1$, and $T$ is its
complementary subspace, $\dim \, T = N -l, \;  T \cap S =
\emptyset$. Let $Y$ be a smooth tangentially nondegenerate and
dually nondegenerate variety of the subspace $T$, $\dim Y =
\mbox{{\rm rank}} \; Y = r < N -l$. Consider an $r$-parameter
family of $l$-dimensional subspaces $L_y = S \wedge y, \; y \in
Y$. This variety  is a \emph{cone} $X$ with vertex $S$ and the
director manifold $Y$. The subspace $T_x (X)$ tangent to the cone
$X$ at a point $x \in L_y (x \notin S)$ is defined by its vertex
$S$ and the subspace $T_y (Y),  T_x (X) = S \wedge T_y (Y)$, and
$T_x (X)$ remains fixed when a point $x$ moves in the subspace
$L_y$. As a result, the cone $X$ is a variety  with a degenerate
Gauss map of dimension $n = l + r$ and rank $r$, with plane
generators $L_y$ of dimension $l$. The generators $L_y$ of the
cone $X$ are leaves of the Monge--Amp\`{e}re foliation associated
with $X$. Note that for a cone $X \subset \mathbb{P}^N$ of rank
$r$ with $l$-dimensional generators, we have $n = r + l, \, n^* =
N - l - 1, \, l^* = n^* - r  = (N - l - 1) - r = N - n - 1$.}

\examp{\label{examp:5} Consider a smooth curve $Y$ in the space
$\mathbb{P}^N$ not belonging to a subspace $\mathbb{P}^{l+1}
\subset \mathbb{P}^N$ and the set of its osculating subspaces
$L_y$ of order and dimension $l$. This set forms a variety $X =
\cup_{y \in Y} L_y$ of dimension $l + 1$ and rank $r = 1$ in
$\mathbb{P}^N$. Such a variety is called a {\em multidimensional
torse}. The subspace $T_y = L_y + \frac{dL_y}{dy}$ is the tangent
subspace to $X$ at all points of its generator $L_y$. Thus, the
subspaces $L_y$ are the leaves of the Monge--Amp\`{e}re foliation
associated with the torse $X$.

Conversely, a   variety  of dimension $n$ and rank 1 is a torse
formed by a family of osculating subspaces of order $n - 1$ of a
curve of class $C^p, \, p \geq n - 1$, in the space
$\mathbb{P}^N$.

Note that for the torse $X$, we have $n = l + 1, \, n^* = N - l -
1,$ and $l^* = N - l - 2$,  i.e., the dual image $X^*$ of a torse
$X$ is a  torse. Note also that the torse considered in Example 1
is a particular case of a multidimensional torse corresponding to
the value $l = 0$.}

\examp{\label{examp:6} Let $Y$ be an $r$-parameter family of
hyperplanes $\xi$ in in the space $\mathbb{P}^{n+1}, \, r < n$,
whose $(n - r)$-dimensional plane generators depend also on $r$
parameters. Then the family $Y$ has an $n$-dimensional envelope
$X$ which is a variety with a degenerate Gauss map  rank $r$. For
this hypersurface, we have   $l = n - r, \, n^* = r, \, l^* = 0$.}

To a cone $X$ of rank $r$ with vertex $S$ of dimension  $l - 1$
(see Example 3), there corresponds a variety  $X^*$ lying in the
subspace $T = S^*, \; \dim T = N - l$.  Because $\dim X^* = n^* =
N - l - 1$, the variety $X^*$ is a hypersurface of rank $r$ in the
subspace $T$. Such a hypersurface was considered in Example 5.

If a tangentially nondegenerate  variety  $X, \, \dim X =
\mbox{{\rm rank}} \, X = r$, belongs to a subspace
$\mathbb{P}^{n+1} \subset \mathbb{P}^N$, then we can consider two
dual maps in the spaces $\mathbb{P}^{n+1}$ and $\mathbb{P}^N$. We
denote the first of these maps by $*$ and the second by $\circ$.
Then under the first map, the image of $X$ is a hypersurface $X^*
\subset \mathbb{P}^{n+1}$, and under the second map, the
hypersurface $X$ is transferred into a cone $X^\circ$ of rank $r$
and dimension $n^\circ = N - n + r - 1$ with an
$(N-n-2)$-dimensional vertex $S = (\mathbb{P}^{n+1})^\circ$ and
$(N-n-1)$-dimensional  plane generators $L^\circ = T (X)^\circ$.

It follows that if the variety $X$ lies in a proper linear
subspace of the space $\mathbb{P}^N$, then Examples 4 and 6 are
mutually dual to each other.

\examp{\label{examp:7} Now we consider the Veronese variety given
as the image of the embedding
$$
V^*: \text{Sym} \; (\mathbb{P}^{2*} \times  \mathbb{P}^{2*})
\rightarrow \mathbb{P}^{5*}
$$
into the projective space $\mathbb{P}^{5*}$. This embedding is
defined by the equations
\begin{equation}\label{eq:8}
x_{ij} = u_i u_j, \;\;\;\; i, j = 0, 1, 2,
\end{equation}
where $u_i$ are projective coordinates in the plane
$\mathbb{P}^{2*}$, i.e., tangential coordinates in the plane
$\mathbb{P}^{2}$, and $x_{ij}$ are projective coordinates in the
space $\mathbb{P}^{5*}$, $x_{ij} = x_{ji}$.

It was shown in [AG 02] (see also [AG 04], p. 77) that the
equation of the variety $V$ that is dual to the variety  $V^*
\subset \mathbb{P}^{5*}$ defined by equations (8) is
\begin{equation}\label{eq:9}
 \det \left(
\begin{array}{lll}x^{00} & x^{01} & x^{02} \\
                  x^{10} & x^{11} & x^{12} \\
                  x^{20} & x^{21} & x^{22}
                  \end{array}\right) = 0, \;\;\;\; x^{ij} =
                  x^{ji}.
\end{equation}

Equation (9) defines in the space $\mathbb{P}^5$ the cubic
hypersurface dual to the Veronese variety (8) and called the {\em
cubic symmetroid}.

The Veronese variety $V^*$ defined by equation (8) is a
tangentially nondegenerate variety in the space $\mathbb{P}^{5*}$.
Thus, by Theorem 1, its dual variety $V$ is a hypersurface with a
degenerate Gauss map of rank two in the space $\mathbb{P}^{5}$
having two-dimensional leaves $L (V)$ of  the Monge--Amp\`{e}re
foliation on $V$. The latter is dual to the tangent bundle $T
(V^*)$ of $V^*$.}

The following table shows the values of dimensions $n$ and $n^*$
of  variety  $X$ and $X^*$, their common rank $r$, and the
dimensions of $l$ and $l^*$ of their plane generators $L$ and
$L^*$ in Examples 1--7:

{\footnotesize
\begin{center}
\hspace*{-15mm}\begin{tabular}{|c||c|c||c|c|c|c|c|c|} \hline
             & $X$  & $N$ & $n$ & $l$ & $r$& $l^*$ & $n^*$ & $X^*$  \\
 \hline \hline
 \textbf{1} & Torse & 3&  $1$ & $0$ & $1$ & $1$  & $2$ & Torse   \\
 \hline
\textbf{2} &Hypersurface&&&&&&&   Tangentially \\
 &  of rank $r$    & $N$ & $N - 1$ & $N - 1 - r$ & $r$ &
 $0$ & $r$ &  nondegenerate \\
&  &&&&&&& variety   \\
\hline
 \textbf{3} & Join & $N \geq 4$&  $3$ & $1$ & $2$ & $N - 4$ & $N - 2$ &    \\
 \hline
  \textbf{4} & Cone & $N$ &  $n$ & $n - r$ & $r$ & $N - n - 1$  & $N - n - 1$ & Hypersurface   \\
 \hline
   \textbf{5} & Multidimen-  &&&&&& & Multidimen-\\
   &sional torse & $N$ &  $l+1$ & $l$ & $1$ & $N - l - 2$  & $N - l - 1$ & sional torse   \\
 \hline
   \textbf{6} & Hypersurface  &&&&&& & \\
   & of rank $r$ & $n + 1$ &  $n$ & $n-r$ & $r$ & $0$  & $r$ & Cone    \\
 \hline
   \textbf{7} & Cubic  &&&&&& & Veronese\\
   & symmetroid & $5$ &  $4$ & $2$ & $2$ & $0$  & $2$ & variety   \\
 \hline
   \end{tabular}
\end{center}
}

\textbf{2.3 Correlative Transformations.} If we have the
identification \linebreak $(\mathbb{P}^N)^* = \mathbb{P}^N$, the
duality principle can be realized by a correlative transformation
of the space $\mathbb{P}^N$.

Consider a  \emph{correlative transformation} ${\cal C}$ (a
\emph{correlation}) in the space $\mathbb{P}^N$ that maps a point
$x \in \mathbb{P}^N$ into a hyperplane $\xi \in \mathbb{P}^N, \;
\xi = {\cal C}(x)$, and preserves the incidence of points and
hyperplanes. A correlation ${\cal C}$ maps a $k$-dimensional
subspace $\mathbb{P}^k \subset \mathbb{P}^N$ into an $(N - k -
1)$-dimensional subspace $\mathbb{P}^{N-k-1} \subset
\mathbb{P}^N$.

We assume that the correlation ${\cal C}$ is nondegenerate, i.e.,
it defines a one-to-one correspondence between points and
hyperplanes of the space $\mathbb{P}^N$.

Analytically, a correlation ${\cal C}$ can be written in the form
$$
\xi_i = c_{ij} x^j, \;\; i, j = 0, 1, \ldots, N,
$$
where $x^i$ are point coordinates and $\xi_i$ are tangential
coordinates in the space $\mathbb{P}^N$. A correlation ${\cal C}$
is nondegenerate if $\det (c_{ij}) \neq 0$.

Consider a smooth curve $C$ in the space $\mathbb{P}^N$ and
suppose that this curve does not belong to a hyperplane. A
correlation ${\cal C}$ maps points of $C$ into hyperplanes forming
a one-parameter family. The hyperplanes of this family envelope a
hypersurface with a degenerate Gauss map of rank one with
$(N-2)$-dimensional generators.

If the curve $C$ lies in a subspace $\mathbb{P}^s \subset
\mathbb{P}^N$, then a  correlation ${\cal C}$ maps points of $C$
into hyperplanes that envelop a hypercone with an
$(N-s-1)$-dimensional vertex.

Further, let $X = V^r$ be an arbitrary tangentially nondegenerate
$r$-dimen\-sional  variety  in the space $\mathbb{P}^N$.  A
correlation ${\cal C}$ maps points of such $V^r$ into hyperplanes
forming an $r$-parameter family. The hyperplanes of this family
envelop a hypersurface $Y = V^{N-1}_r$ with a degenerate Gauss map
of rank $r$. The  generators of this hypersurface $X$ are of
dimension $N-r-1$ and correspond to the tangent subspaces $T_x
(V^r)$.

If the  tangentially nondegenerate variety  $V^r$ belongs to a
subspace $\mathbb{P}^s \subset \nolinebreak \mathbb{P}^N,
\linebreak s > r$, then  the hypersurface $Y = V^{N-1}_r$
corresponding to $V^r$ under a correlation ${\cal C}$ is a
hypercone with an $(N-s-1)$-dimensional vertex.

Now let $X = V^n_r$ be a   variety  with a degenerate Gauss map of
rank $r$. Then we can prove the following result, which fully
corresponds to Theorem 1.

\setcounter{theorem}{1}
\begin{theorem}
A  correlation ${\cal C}$ maps an $n$-dimensional dually
nondegenerate variety $X = V^n_r$ with a degenerate Gauss map of
rank $r$ with plane generators of dimension $l = n - r$ into an
$(N-l-1)$-dimensional variety $X^* = V^{N-l-1}_r$, with a
degenerate Gauss map of the same rank $r$ with
$(N-n-1)$-dimensional plane generators.
\end{theorem}

\begin{proof}
A correlation ${\cal C}$ sends an $l$-dimensional plane generator
$L \subset X$ to an $(N-l-1)$-dimensional plane
$\mathbb{P}^{N-l-1}$, and a tangent subspace $T_x (X)$ to an
$(N-n-1)$-dimensional plane  $\mathbb{P}^{N-n-1}$, where
$\mathbb{P}^{N-n-1} \subset  \mathbb{P}^{N-l-1}$. Because  both of
these planes depend on $r$ parameters, the planes
$\mathbb{P}^{N-n-1}$ are generators of the  variety  ${\cal C}
(X)$, and the planes $\mathbb{P}^{N-l-1}$ are its tangent
subspaces. Thus, the variety ${\cal C} (X)$ is a  variety $X^* =
V^{N-l-1}_r$ of dimension $N - l - 1$ and rank $r$.
\end{proof}

\section{Basic Equations of a Variety with  \\ a Degenerate Gauss Map.}

In this section, we find the basic equations of a variety $X$ with
a degenerate Gauss map of dimension $n$ and rank $r$ in a
projective space $\mathbb{P}^N$.

In what follows, we will use the following ranges of indices:
$$
a, b, c = 1, \ldots, l; \; p, q = l + 1, \ldots , n; \; \alpha,
\beta = n + 1, \ldots , N.
$$

A point $x \in X$ is said to be a \textit{regular point} of the
map $\gamma$ and the variety $X$ if $\dim T_x X = \dim X = n$, and
a point $x \in X$ is called a \textit{singular point} of the leaf
$L \subset X$ if $\dim T_x X
> \dim X = n$.

In what follows, we assume that every plane generator $L$ of a
variety $X$ with a degenerate Gauss map has at least one regular
point. Otherwise (i.e., if all points of $L$ are singular) the
Monge--Amp\`{e}re foliation is degenerate, and we will not
consider this case.

We associate a family of moving frames $\{A_u\}, \, u = 0, 1,
\ldots , N$, with $X$ in such a way that the point $A_0 = x$ is a
regular point of $X$; the points $A_a$ belong to the leaf $L$ of
the Monge--Amp\`{e}re foliation passing through the point $A_0$;
the points $A_p$ together with the points $A_0, A_a$ define the
tangent subspace $T_L X$ to $X$; and the points $A_\alpha$ are
located outside the subspace $T_L X$.

The equations of infinitesimal displacement of the moving frame
$\{A_u\}$ are
 \begin{equation}\label{eq:10}
 dA_u = \omega_u^v A_v, \;\;\;\;u, v = 0, 1, \ldots , N,
\end{equation}
where $\omega_u^v$ are 1-forms satisfying the structure equations
of the projective space $\mathbb{P}^N$:
 \begin{equation}\label{eq:11}
 d \omega_u^v = \omega_u^w \wedge \omega_w^v, \;\;\;\;u, v, w = 0, 1, \ldots , N.
\end{equation}

As a result of the specialization of the moving frame mentioned
above, we obtain the following equations of the variety $X$ (see
[AG 04], Section 3.1):
\begin{equation}\label{eq:12}
\omega_0^\alpha = 0,
\end{equation}
\begin{equation}\label{eq:13}
\omega_a^\alpha = 0,
\end{equation}
\begin{equation}\label{eq:14}
\omega_p^\alpha = b_{pq}^\alpha \omega^q, \;\; b_{pq}^\alpha =
b_{qp}^\alpha,
\end{equation}
\begin{equation}\label{eq:15}
\omega^p_a = c^p_{aq} \omega^q.
\end{equation}

The 1-forms $\omega^q := \omega_0^q$ in these equations are basis
forms of the Gauss image $\gamma (X)$ of the variety $X$, and the
quantities $b_{pq}^\alpha$ form the second fundamental tensor of
the variety $X$ at the point $x$. The quantities $b_{pq}^\alpha$
and $c^p_{aq}$ are related by the following equations:
\begin{equation}\label{eq:16}
b_{sq}^\alpha c^s_{ap} = b_{sp}^\alpha c^s_{aq}.
\end{equation}

Equations (14) and (15) are called the \emph{basic equations} of a
variety $X$ with a degenerate Gauss map (see [AG 04], Section
3.1).

Note that  under transformations of the points $A_p$, the
quantities $c^p_{aq}$ are transformed as tensors. As to the index
$a$, the  quantities $c^p_{aq}$ do not form a tensor with respect
to this index. Nevertheless, under transformations of the points
$A_0$ and $A_a$, the quantities $c^p_{aq}$ along with the unit
tensor $\delta_q^p$ are  transformed as  tensors. For this reason,
the system of quantities $c^p_{aq}$ is called a {\em quasitensor}.

By (12) and (13), the equations of infinitesimal displacement of
the moving frame associated with a variety $X$ with a degenerate
Gauss map have the form
\begin{equation}\label{eq:17}
\renewcommand{\arraystretch}{1.5}
\left\{
\begin{array}{ll}
dA_0 = \omega_0^0 A_0 + \omega^a A_a + \omega^p A_p, \\
dA_a = \omega_a^0 A_0 + \omega^b_a A_b + \omega^p_a A_p, \\
dA_p = \omega_p^0 A_0 + \omega^a_p A_a + \omega^q_p A_p
       + \omega^\alpha_p A_\alpha, \\
dA_\alpha = \omega_\alpha^0 A_0 + \omega^\alpha_a A_a +
\omega^\alpha_q A_q + \omega_\alpha^\beta A_\beta,
\end{array}
\right.
\renewcommand{\arraystretch}{1}
\end{equation}
where here and in what follows, unless otherwise stated, the
indices take the values indicated earlier.

Denote by $B^\alpha$ and $C_a$  the $(r \times r)$-matrices of
coefficients occurring in equations (14) and (15):
$$
B^\alpha = (b^\alpha_{pq}), \;\; C_a = ( c^p_{aq}).
$$
Sometimes we will use the identity matrix $C_0 = (\delta^p_q)$ and
the index \linebreak $i = 0, 1, \ldots, l$, i.e., $\{i\} = \{0,
a\}$. Then equations (14) and (16) can be combined and written as
follows:
$$
(B^\alpha C_i)^T = (B^\alpha C_i),
$$
i.e., the matrices
$$
H^\alpha_i = B^\alpha C_i = (b^\alpha_{qs} c^s_{ip})
$$
are symmetric.

The \textit{second fundamental form} $II$ of the variety $X$ at
the point $x$ is the linear span of the system of quadratic forms
\begin{equation}\label{eq:18}
\Phi^\alpha = b_{pq}^\alpha \omega^p \omega^q
\end{equation}
(see, for example, [AG 04], Section 2.2; [FP 01], Section 2.4; [GH
79]; [L 94]; and [L 99], Section 4.1.5).

Suppose now that $x = x^i A_i, \, i = 0, 1, \dots, l$, is an
arbitrary point of a generator $L$ of the variety $X$. By (17),
the differential of this point has the form
$$
d x = (d x^i + x^j \omega_j^i) A_i + x^i c^p_{iq} \omega^q A_p,
$$
where $i, j = 0, 1, \dots, l$ and $c^p_{0q} = \delta^p_q$. The
matrix
$$
J^p_q = (x^i c^p_{iq})
$$
is called the \textit{Jacobi matrix} of the variety $X$. A point
$x \in X$ is \textit{regular} if $\det (J^p_q) \neq 0$ at $x$, and
a point $x \in X$ is \textit{singular} if $\det (J^p_q) = 0$ at
$x$, i.e.,
\begin{equation}\label{eq:19}
\det (x^i c^p_{iq}) = 0
\end{equation}
 (see [AG 04, Section 3.2). The set of
singular points forms the focus algebraic hypersurface of order
$r$ on $L$.

Consider the bundle $\Xi$ of tangent hyperplanes $\xi$ of the
variety $X$ passing through its tangent subspace $T_L$. Because
$\dim T_L = n$, the dimension of the bundle $\Xi$ is $N - n - 1,
\, \dim \Xi = N - n - 1$. In the space $\mathbb{P}^N$, a tangent
hyperplane $\xi$ is determined by the equation $\xi = \xi_\alpha
x^\alpha = 0$. In the linear system of quadratic forms defined by
the forms $\Phi^\alpha = b^\alpha_{pq} \omega^p \omega^q$, to the
tangent hyperplane $\xi = (\xi_\alpha)$, there corresponds the
quadratic form
$$
\Phi (\xi) = \xi_\alpha b^\alpha_{pq} \omega^p \omega^q.
$$
The tangent hyperplane $\xi$ is called \textit{singular} if the
quadratic form $\Phi (\xi)$ is singular, i.e., if
\begin{equation}\label{eq:20}
\det (\xi_\alpha b^\alpha_{pq}) = 0.
\end{equation}

\section{The Generalization of a Theorem of Griffiths--Harris}

\textbf{3.1 A  Theorem of Griffiths--Harris.} Griffiths and Harris
in [GH 79] (see Theorem 3.5 in [GH 79]) proved the following
theorem:

\begin{theorem}
If a variety $X$ is tangentially nondegenerate, then its dual
variety $X^*$ is dually degenerate if and only if at any smooth
point $x \in X$ every quadratic form $\Phi = \xi_\alpha
\Phi^\alpha$  belonging to the  second fundamental form $II$ of
the variety $X$ is singular.
\end{theorem}

The following example illustrates this theorem.

\setcounter{theorem}{7}

\examp{\label{examp:8} \textit{The Segre variety}   $S (m, n)$ is
the embedding of the direct product of the projective spaces
$\mathbb{P}^m$ and $\mathbb{P}^n$ in the space $\mathbb{P}^{mn + m
+ n}$ (see  [GH 79] or  [T 01]):
$$
S: \mathbb{P}^m \times \mathbb{P}^n \rightarrow \mathbb{P}^{mn + m
+ n},
$$
defined by the equations $$z^{ik} = x^i y^k,$$ where $i = 0, 1,
\ldots, m, \, k = 0, 1, \ldots, n,$ and $x^i, y^k,$ and $z^{ik}$
are the coordinates of points in the spaces $\mathbb{P}^m,
\mathbb{P}^n,$ and $\mathbb{P}^{mn + m + n}$, respectively. This
manifold has the dimension $m + n,\, \dim S (m, n) = m + n$.

Consider in the spaces $\mathbb{P}^m$ and $\mathbb{P}^n$
projective frames $\{A_0, A_1, \ldots, A_m\}$ and $\{B_0, B_1,
\ldots, B_n\}$. Then in the space $\mathbb{P}^{mn+m+n}$ we obtain
the projective frame
$$
\{A_0 \otimes B_0, \, A_0 \otimes B_k, \, A_i \otimes B_0, \, A_i
\otimes B_k\}
$$
(here and in what follows $i, j = 1, \ldots, m; \, k, l = 1,
\ldots , n$) consisting of \linebreak $(m + 1)(n+1)$ linearly
independent points of the space $\mathbb{P}^{mn+m+n}$. The point
$A_0 \otimes B_0$ is the generic point of the variety $S$.

In the  spaces $\mathbb{P}^m$ and $\mathbb{P}^n$, we have the
following equations of infinitesimal displacements of the moving
frames $\{A_0, A_1, \ldots, A_m\}$ and $\{B_0, B_1, \ldots,
B_n\}$:
$$
d A_0 = \omega_0^0 A_0 + \omega_0^i A_i, \;\; d B_0 = \sigma_0^0
B_0 + \sigma_0^k B_k.
$$
Hence
$$
d (A_0 \otimes B_0) = (\omega_0^0 + \sigma_0^0) (A_0 \otimes B_0)
+ \omega_0^i (A_i \otimes B_0) + \sigma^k_0  (A_0 \otimes B_k),
$$
and the subspace in $\mathbb{P}^{mn+m+n}$ spanned by the points
$A_0 \otimes B_0, A_i \otimes B_0$, and $A_0 \otimes B_k$ is the
tangent subspace to the Segre variety $S$ at the point $A_0
\otimes B_0$:
$$
T_{A_0 \otimes B_0} = \text{Span} \; (A_0 \otimes B_0, \, A_i
\otimes B_0, \, A_0 \otimes B_k).
$$

The second differential of the point $A_0 \otimes B_0$ has the
form:
$$
d^2 (A_0 \otimes B_0) = 2 \, \omega_0^i \, \sigma_0^k \, A_i
\otimes B_k \pmod{T_{A_0 \otimes B_0}}.
$$
Hence the osculating subspace $T^2_{{A_0 \otimes B_0}} (S)$ to the
variety $S$ coincides with the entire space
$\mathbb{P}^{mn+m+n}/T_{A_0 \otimes B_0}$, and its second
fundamental form is the linear span of the system of bilinear
forms
$$
\Phi^{ik} = \omega_0^i \,  \sigma_0^k.
$$
The total number of these forms is $mn$. The equations $\omega_0^i
= 0$ determine \linebreak $n$-dimensional plane generators on $S$,
and the equations $\sigma_0^k = 0$ determine its $m$-dimensional
plane generators.

Consider a tangent hyperplane to the Segre variety $S$ at the
point $A_0 \otimes B_0$. Because such a hyperplane contains the
tangent subspace $T_{A_0 \otimes B_0}$, its equation can be
written in the form
$$
\xi = \xi_{ik} z^{ik} = 0,
$$
where $i = 1, \ldots, m; \, k = 1, \ldots , n$, and $z^{ik}$ are
coordinates of points in the space $\mathbb{P}^{mn+m+n}/T_{A_0
\otimes B_0} $. As a result, the second fundamental form of the
variety $S$ with respect to the hyperplane $\xi$ is
$$
\Phi (\xi) = \xi_{ik}\, \omega_0^i \, \sigma_0^k
$$
(cf. equations (2.21) in [AG 04]). The forms $\Phi (\xi)$ is a
linear combination of the linearly independent forms $\Phi^{ik}$,
and the matrix of this form is
$$
\Xi = \displaystyle \frac{1}{2} \left(
\renewcommand{\arraystretch}{1.3}
\begin{array}{ll}
0 & (\xi_{ik}) \\
(\xi_{ki}) & 0
\end{array}
\renewcommand{\arraystretch}{1}
\right).
$$
In this formula $(\xi_{ik})$ is a rectangular $(m\times n)$-matrix
and  $(\xi_{ki})$ is its transpose.

It follows that $\det \Xi = 0$ if $m \neq n$. In this case, all
quadratic forms belonging to the second fundamental forms of the
variety $S$ are singular, and the dual defect $\delta_* (S)$ of
$S$ equals $| n - m |: \, \delta_* (S) = | n - m |$. If $m \neq
n$, then the variety $S$ is dually degenerate. The variety $S$ is
dually nondegenerate if and only if $m = n$. }

For other proofs of the formula $\delta_* (S) = | n - m |$ see [L
99], p. 52; or [FP 01], p. 110; or [Ha 92], p. 198; or [AG 04], p.
74--76).

\vspace*{2mm}

\textbf{3.2 The Generalization of a Theorem of Griffiths--Harris.}
In Section 1.2, we defined the dual variety $X^* \subset
(\mathbb{P}^{N})^*$ for a variety $X \subset \mathbb{P}^N$ with a
degenerate Gauss map of dimension $n$ and rank $r$ as the set of
tangent hyperplanes \linebreak $\xi \, (\xi \supset T_L X)$ to
$X$. It follows that the dual variety $X^*$ is a fibration whose
fiber is the bundle
$$
\Xi = \{\xi | \xi \supset T_L X\}
$$
of hyperplanes $\xi$ containing the tangent subspace $T_L X$ and
whose base is the manifold
$$
B = X^* / \Xi.
$$
As we noted in Section 1.2, the dimension of a fiber $\Xi$ of this
fibration equals $N - n - 1$, and the dimension of the base $B$
equals $r, \dim B = \nolinebreak r$, i.e., the dimension of $B$
coincides with the rank of the variety $X$. This implies that in
the general case,
$$
\dim X^* = (N - n - 1) + r = N - l - 1
$$
(cf. formula (6)).

As was noted in Section 2, for a dually degenerate variety $X$
with a degenerate Gauss map, we have
$$
\dim X^* < N - l - 1.
$$

The following theorem generalizing Theorem 3 expresses this
condition in terms of the second fundamental forms of the variety
$X$.

\setcounter{theorem}{3}
\begin{theorem}
The dual variety $X^* \subset (\mathbb{P}^{N})^*$ of a variety $X$
with a degenerate Gauss map is dually degenerate if and only if at
any smooth point $x \in X$ every quadratic form $\Phi = \xi_\alpha
\Phi^\alpha$  belonging to the  second fundamental form $II$ of
the variety $X$ is singular.
\end{theorem}

\begin{proof} Consider the bundle  ${\cal R}
(X)$ of frames associated with a variety $X$ with a degenerate
Gauss map, which we constructed earlier in this section. The basis
forms of the bundle ${\cal R}(X)$, as well as the basis forms of
the tangent bundle $T (X)$ and the Monge--Amp\`{e}re foliation of
the variety $X$, are also called the \textit{horizontal forms},
and the secondary forms of all these bundles are called the
\textit{fiber} or \textit{vertical forms}. The horizontal forms
$\omega^p, \, p = l + 1, \dots, n$, are linearly independent, and
their number equals $r$. Thus, these forms are basis forms in the
bundle ${\cal R}(X)$.
 On the bundle ${\cal R} (X)$ the equations of
infinitesimal displacement of a frame have the form (17).

In this proof we will use the following ranges of indices:
$$
\renewcommand{\arraystretch}{1.3}
\begin{array}{ll}
0 \leq u, v \leq N, & 1 \leq i, j \leq n, \\ 1 \leq a, b \leq l, &
l+1 \leq p, q \leq n,  \\  n+1 \leq \alpha, \beta \leq N, &  n+1
\leq \rho, \sigma \leq N-1.
\end{array}
\renewcommand{\arraystretch}{1}
$$

Consider now  the dual coframe (or tangential frame)
$\{\alpha^u\}$ in the space $(\mathbb{P}^{N})^*$  to the frame
$\{A_u\}$ (see Section 1.3 in [AG 04]). The hyperplanes $\alpha^u$
of the frame  $\{\alpha^u\}$  are connected with the points of the
frame $\{ A_u\}$ by the conditions
\begin{equation}\label{eq:21}
(\alpha^u, A_v) = \delta^u_v.
\end{equation}
Conditions (21) mean that the hyperplane $\alpha^u$ contains all
points $A_v, \; v \neq u$, and that the condition of normalization
$(\alpha^u, A_u) = 1$ holds.

The equations of infinitesimal displacement of the tangential
frame $\{\alpha^u\}$ have the form (see Section 1.3 in [AG 04])
\begin{equation}\label{eq:22}
d \alpha^u = \widetilde{\omega}^u_v \alpha^v, \;\; u, v = 0, 1,
\dots , N,
\end{equation}
where the forms $\widetilde{\omega}_u^v$ are related to the forms
$\omega_u^v$ by the following formulas:
$$
\widetilde{\omega}^u_v = - \omega^u_v.
$$
Hence equations (22) can be written as
\begin{equation}\label{eq:23}
d \alpha^u = - \omega^u_v \alpha^v.
\end{equation}

Recalling that
$$
\renewcommand{\arraystretch}{1.3}
\left\{
\begin{array}{ll}
dA_0 \equiv \omega^a A_a + \omega^p A_p \pmod{A_0}, \\
dA_a \equiv \omega_a^p A_p \pmod{A_0, A_1, \ldots, A_l}, \\
dA_p \equiv \omega_p^\sigma A_\sigma + \omega_p^N A_N \pmod{A_0,
A_1,  \ldots, A_n}
\end{array} \right.
\renewcommand{\arraystretch}{1}
$$
(cf. equations (17)) and
$$
\omega^\alpha = 0, \;\; \omega_a^\alpha = 0, \;\; \omega_p^\alpha
= b_{pq}^\alpha \omega^p
$$
(cf. equations (12), (13), and (14)), we find from (23) that
$$
d\alpha^N \equiv - \omega^N_a \alpha^a -  \omega_p^N \alpha^p -
\omega_\sigma^N \alpha^{\sigma} \pmod{\alpha^N}.
$$
The $N-n-1$ forms $\omega_\sigma^N$ determine the infinitesimal
displacement of the hyperplane $\xi = \alpha^N$ in the bundle
$\Xi$ of tangent hyperplanes $\xi$ containing the tangent subspace
$T_L X$, i.e., these forms are the fiber forms on the dual variety
$X^*$. The number $N - n - 1$ coincides with the dimension of a
fiber of this bundle. Hence, forms $\omega_\sigma^N$ are linearly
independent.

A basis of the fibration $X^*$ is the span $S^N$ of the forms
$\omega_p^N$, i.e., the forms $\omega_p^N$ are horizontal on
$X^*$. Because
$$
\omega_p^N = b_{pq}^N \omega^q, \;\;  b_{pq}^N = b_{qp}^N, \;\; p,
q = l + 1, \ldots , n,
$$
the dimension of $S^N$ does not exceed the rank $r = n - l$ of the
variety $X$.

Consider the exterior product
$$
\omega_{l+1}^N \wedge \ldots \wedge \omega_n^N = \det \,(b_{pq}^N)
\, \omega^{l+1} \wedge \ldots \wedge \omega^n.
$$
It is easy to see that $\dim \, S^N = \text{rank} \; (b^N_{pq})$,
and $\dim \, S^N < r$ if and only if $\det \,(b_{pq}^N) = 0$.

Because $\alpha^N$ was any of the hyperplanes $\alpha^\beta$, we
have
$$
\det \,(b_{pq}^\beta) = 0.
$$

Moreover, the tangent hyperplane $\xi$ can be chosen arbitrarily
from the system $\xi = \xi_\beta \alpha^\beta$. This system of
tangent hyperplanes passing through the tangent subspace $T_L X$
determines the system of quadratic forms
\begin{equation}\label{eq:24}
\xi_\beta b_{pq}^\beta \omega^p \omega^q
\end{equation}
whose span is the second fundamental form $II$ of the variety $X$
and the system of second fundamental tensors
$$
 \xi_\beta b_{pq}^\beta
$$
of this variety $X$. This proves the theorem statement: the
variety $X$ is dually degenerate if and only if every quadratic
form $\Phi = \xi_\alpha \Phi^\alpha$  belonging to the  second
fundamental form $II$ of the variety $X$ is singular.
\end{proof}

Note that if $r = n$ $($i.e., if a variety $X$ is tangentially
nondegenerate$)$, then  we obtain Theorem 3.

We emphasize that unlike in (24), the basis forms in Theorem 3 are
the forms $\omega^i, \, i = 1, \dots, n$.

\begin{corollary}
A variety $X$ with a degenerate Gauss map is dually nondegenerate
$($i.e., the dimension of its dual variety $X^* \subset
(\mathbb{P}^{N})^*$ equals $N - l - 1)$ if and only if at any
smooth point $x \in X$ there is at least one nonsingular quadratic
form $\Phi = \xi_\alpha \Phi^\alpha$ belonging to the second
fundamental form $II$ of the variety $X$.
\end{corollary}

\setcounter{theorem}{8}

\examp{\label{examp:9} It is not difficult to construct an example
of a dually degenerate variety $X$ with a degenerate Gauss map
departing from Example 8 of a Segre variety $S (m, n) =
\mathbb{P}^m \times \mathbb{P}^n \subset \mathbb{P}^{mn + m + n}$.
To this end, we consider the projective space $\mathbb{P}^N$ of
dimension $N = mn + m + n + l$ and its two complementary subspaces
$\mathbb{P}^{mn + m + n}$ and $\mathbb{P}^{l-1}$. Consider the
Segre variety $S (m, n)$ in the first of these subspaces
$\mathbb{P}^{mn + m + n}$, and let $C (m, n)$ be a cone with
vertex $\mathbb{P}^{l-1}$ whose director variety is the variety $S
(m, n)$. The cone  $C (m, n)$ has a degenerate Gauss map of rank
$r = m + n = \dim S (m, n)$. The quadratic forms belonging to the
second fundamental form  of the cone  $C (m, n)$ coincide with
quadratic forms belonging to the second fundamental form of its
director variety $S (m, n)$. Because the variety $S (m, n)$ is
dually degenerate if $m \neq n$, the cone  $C (m, n)$ is also
dually degenerate. }

\noindent {\em Authors' addresses}:\\

\noindent
\begin{tabular}{ll}
M.~A. Akivis &V.~V. Goldberg\\
Department of Mathematics &Department of Mathematical Sciences\\
Jerusalem College of Technology---Mahon Lev &  New
Jersey Institute of Technology \\
Havaad Haleumi St., P. O. B. 16031 & University Heights \\
Jerusalem 91160, Israel &  Newark, N.J. 07102, U.S.A. \\
 & \\
 E-mail address: akivis@mail.jct.ac.il & E-mail address:
 vlgold@oak.njit.edu
 \end{tabular}

\end{document}